\documentclass[10 pt]{amsart}

\usepackage[colorlinks=true,linkcolor=red,citecolor=blue,urlcolor=blue,hypertexnames=false]{hyperref}
\usepackage{bookmark}
\usepackage{amsthm,thmtools,amssymb,amsmath,amscd}

\usepackage{fancyhdr}
\usepackage{esint}
\usepackage{enumerate}

\usepackage{pictexwd,dcpic}
\swapnumbers
\declaretheorem[name=Theorem,numberwithin=section]{thm}
\declaretheorem[name=Remark,style=remark,sibling=thm]{rem}
\declaretheorem[name=Lemma,sibling=thm]{lemma}
\declaretheorem[name=Proposition,sibling=thm]{prop}
\declaretheorem[name=Definition,style=definition,sibling=thm]{defn}
\declaretheorem[name=Corollary,sibling=thm]{cor}

\numberwithin{equation}{section}

\usepackage{cleveref}
\crefname{lemma}{Lemma}{Lemmata}
\crefname{prop}{Proposition}{Propositions}
\crefname{thm}{Theorem}{Theorems}
\crefname{cor}{Corollary}{Corollaries}
\crefname{defn}{Definition}{Definitions}
\crefname{example}{Example}{Examples}
\crefname{rem}{Remark}{Remarks}
\crefname{assum}{Assumption}{Assumptions}
\crefname{notation}{Notation}{Notation}

\usepackage{autonum}

\newcommand{\ti}{\tilde}

\newcommand{\wh}{\widehat}
\newcommand{\bs}{\backslash}
\newcommand{\cn}{\colon}
\newcommand{\sub}{\subset}

\newcommand{\R}{\mathbb{R}}

\newcommand{\bbS}{\mathbb{S}}
\newcommand{\bbR}{\mathbb{R}}

\newcommand{\B}{\mathbb{B}}
\newcommand{\8}{\infty}

\newcommand{\al}{\alpha}

\newcommand{\e}{\epsilon}
\newcommand{\ka}{\kappa}
\newcommand{\la}{\lambda}

\newcommand{\om}{\omega}
\newcommand{\s}{\sigma}
\newcommand{\Si}{\Sigma}
\newcommand{\p}{\varphi}

\newcommand{\cL}{\mathcal{L}}
\newcommand{\cT}{\mathcal{T}}

\newcommand{\cW}{\mathcal{W}}
\newcommand{\cH}{\mathcal{H}}

\newcommand{\del}{\partial}
\newcommand{\n}{\nabla}

\newcommand{\fa}{\forall}

\newcommand{\ip}[2]{\left\langle #1,#2 \right\rangle}
\newcommand{\fr}[2]{\frac{#1}{#2}}
\newcommand{\x}{\times}

\DeclareMathOperator{\const}{const}

\DeclareMathOperator{\Rm}{Rm}

\newcommand{\pf}[1]{\begin{proof} #1 \end{proof}}
\newcommand{\eq}[1]{\begin{equation}\begin{alignedat}{2} #1 \end{alignedat}\end{equation}}

\newcommand{\br}[1]{\left(#1\right)}
\newcommand{\enum}[1]{\begin{enumerate}[(i)] #1 \end{enumerate}}

\newcommand{\ra}{\rightarrow}
\newcommand{\hra}{\hookrightarrow}

\newcommand{\hp}{\hphantom}
\newcommand{\q}{\quad}

\begin{document}

\title[Alexandrov-Fenchel inequalities with free boundary]{Alexandrov-Fenchel inequalities for convex hypersurfaces with free boundary in a ball}
\author[J. Scheuer]{Julian Scheuer}
\address{Albert-Ludwigs-Universit\"{a}t,
Mathematisches Institut, Abteilung Reine Mathematik, Ernst-Zermelo-Str. 1, 79104
Freiburg, Germany}
\email{julian.scheuer@math.uni-freiburg.de}
\author[G. Wang]{Guofang Wang}
\address{Albert-Ludwigs-Universit\"{a}t,
Mathematisches Institut, Abteilung Reine Mathematik, Ernst-Zermelo-Str. 1, 79104
Freiburg, Germany}
\email{guofang.wang@math.uni-freiburg.de}
\author[C. Xia]{Chao Xia}
\address{School of Mathematical Sciences, Xiamen University, 361005,
Xiamen, P.R. China}
\email{chaoxia@xmu.edu.cn}
\date{\today}

\begin{abstract} In this paper we first introduce  quermassintegrals for free boundary hypersurfaces in the $(n+1)$-dimensional Euclidean unit ball. 
Then we solve some related isoperimetric type problems for convex free boundary hypersurfaces, which lead 
to new Alexandrov-Fenchel inequalities.
In particular, for $n=2$ we obtain a Minkowski-type inequality and for $n=3$ we obtain an optimal Willmore-type inequality. 
To prove these estimates, we employ a specifically designed locally constrained inverse harmonic mean curvature flow with free boundary.
\end{abstract}

\maketitle

\section{Introduction}

The study of free boundary surfaces or hypersurfaces has a very long history.
It goes back at least to  Courant \cite{Courant:12/1940}. Since then, there have been a lot research activities on this topic,
see for example \cite{GruterJost:/1986,GruterJost:09/1986,Smyth:/1984,Struwe:/1984}. 
Due to recent developments inspired by the proof of the Willmore conjecture by Marques-Neves \cite{MarquesNeves:03/2014}
and Fraser-Schoen's work on the first Steklov eigenvalue and minimal 
free boundary surfaces \cite{FraserSchoen:03/2011, FraserSchoen:03/2016}, there are many new interesting results on the existence of minimal free boundary surfaces. We are interested in the study of hypersurfaces in the unit ball with free boundary in the unit sphere.
Such hypersurfaces share many properties with closed ones in the Euclidean space or the sphere. This naturally leads to questions, such as if a certain result is true for hypersurfaces of the sphere, is it also true for hypersurfaces with free boundary on the sphere. One good example is a result of
Fraser-Schoen \cite{FraserSchoen:03/2016}: 

\begin{center}{\it 
If $\Sigma$ is a free boundary minimal surface given by first Steklov  eigenfunctions, homeomorphic to the annulus, then $\Si$  is a critical catenoid.}
\end{center}

This is a free boundary version of a result of Montiel-Ros \cite{MontielRos:/1986} for closed minimal surfaces in the unit sphere. 
However, the proof is much more delicate. In fact, in many cases the free boundary version is more difficult or even left open.
For example, a free  boundary version of the Lawson conjecture solved by Brendle \cite{Brendle:/2013b} is still open \cite{FraserLi:/2014,Nitsche:03/1985}:
Is the critical catenoid the unique embedded annulus? One can also ask if  the critical catenoid  has the least area among minimal annuli.
Its counterpart in the unit sphere was solved by Marques-Neves \cite{MarquesNeves:03/2014} and plays a key role in the resolution of the Willmore conjecture.
We remark  that the critical catenoid does not minimise the Willmore functional 
(a conformally invariant functional with a correction term involving the geodesic curvature of the boundary), but the piece of  the Clifford torus 
intersecting
the unit sphere does. The difficulties arise from two different geometries. Unlike the theory of closed hypersurfaces, a free boundary hypersurface
lies in two different geometries: the Euclidean geo\-metry for its interior and the spherical geometry for its boundary. 
One might have to deal with these two geometries simultaneously. In our recent work, we have used a vector field which is conformal Killing 
in the Euclidean space
and is also conformal Killing on the unit sphere  after restriction. This   vector field  helped us to solve the stability problem for CMC free boundary hypersurfaces
in \cite{WangXia:08/2017}. It will also play an important role in this paper. The objective of this paper is to find a counterpart of quermassintegrals 
(cross-section integrals)
for free boundary hypersurfaces and consider related isoperimetric type problems. We will see in this paper  that 
the quermassintegrals  we find involve not only geometric integrals in the Euclidean space, but also geometric integrals in the spherical space.
We hope that our research provides a deeper understanding of free boundary hypersurfaces and leads to more interesting questions, for example:
\begin{center}{\it{Is there
an integral geometry for free boundary hypersurfaces in a Euclidean ball?}}
\end{center}

This article is mainly about convex hypersurfaces with boundary $(\Si,\del\Si)$ in the closed $(n+1)$-dimensional Euclidean unit ball $\bar\B^{n+1}$. 
Hereby we mean that the second fundamental form is semi-definite. At the boundary $\del\Si$ the hypersurface is supposed to meet the unit sphere $\bbS^{n}\sub\R^{n+1}$ perpendicularly. We will label this property via the following terminology.

\begin{defn}{\label{Free boundary}}
Let $n\geq 2$ and $\Si\sub\bar\B^{n+1}$ be a smooth compact, embedded topological closed $n$-ball, given by an immersion 
\eq{x\cn \bar{\B}^{n}&\ra \Si\sub \bar\B^{n+1}.} We say that $\Si$ has {\it{free boundary in the unit ball}}, if
\eq{x(\B^{n})&\sub \B^{n+1}\\
	\del\Si&\sub\bbS^{n}\\
		\ip{\bar N}{\nu}_{|\del \B^{n}}&=0,}
where $\nu$ is a smooth choice of a unit normal field on $\Si$ and $\bar N$ is the position vector field in $\R^{n+1}$, i.e. its restriction to $\bbS^{n}$ is the outward normal on $\bbS^{n}$. We also identify, without ambiguity, the outward pointing conormal $\mu$ of $\del\Si\sub\Si$ with $\bar N$.
\end{defn}

We are interested in geometric inequalities for free boundary hypersurfaces in the unit ball. The most prominent and best studied one is the relative isoperimetric inequality, which holds in a broad sense, cf. \cite[Thm.~18.1.3]{BuragoZalgaller:/1988}. Here we want to prove higher order generalisations of this, which are classically called Alexandrov-Fenchel inequalities. In the case of closed convex hypersurfaces of the Euclidean space $\R^{n+1}$ these are
\eq{\label{Alex}\int_{\Si}H_{k}\geq \om_{n}^{\fr{k-l}{n-l}}\br{\int_{\Si}H_{l}}^{\fr{n-k}{n-l}},\q 0\leq l< k\leq n,}
with equality precisely at spheres. $\om_{n}$ denotes the surface area of $\bbS^{n}$ and $H_{k}$ the normalised $k$-th mean curvatures of $\Si$ with the convention $H_{0}=1$ and $H_{n+1}=0$, see \cite{Schneider:/2014}. There are generalisations to other space forms, such as the hyperbolic space \cite{GeWangWu:04/2014,GeWangWu:/2013,LiWeiXiong:03/2014,WangXia:07/2014} and the sphere \cite{MakowskiScheuer:11/2016,WeiXiong:/2015}. The aim of the present paper is to find the correct extensions to hypersurfaces with free boundary in the unit ball.

In order to state the main theorem, we have to introduce some more notation.
If $\Si$ is strictly convex, then $\del\Si\sub\bbS^{n}$ is a strictly convex hypersurface of the unit sphere and bounds a 
strictly convex body in $\bbS^{n}$, which we denote by $\wh{\del\Si}$, cf. \cite{CarmoWarner:/1970}. 
Denote by $\wh\Si$ the domain enclosed by $\Si$ in $\B^{n+1}$, which contains $\wh{\del \Si}$. 
Let $\s_{k}$ denote the $k$-th elementary symmetric polynomial, evaluated at the principal curvatures of the hypersurface 
$\Si$,
\eq{\s_{k}=\binom{n}{k}H_{k},\q 0\leq k\leq n, \q \s_{n+1}=H_{n+1}=0.}
Then we define the following geometric functionals, which we expect to be the correct counterparts to the quermassintegrals for closed convex hypersurfaces.

\eq{\label{Geom-Quant}W_{0}(\wh\Si)&=|\wh\Si|\\
	W_{1}(\wh\Si)&=\fr{1}{n+1}|\Si|\\
	W_{k}(\wh\Si)&=\fr{1}{n+1}\int_{\Si}H_{k-1}+\fr{k-1}{(n+1)(n-k+2)}W_{k-2}^{\bbS^{n}}(\wh{\del\Si}),\q 2\leq k\leq n+1,}
where for a $k$-dimensional submanifold $M\sub\R^{n+1}$ (with or without boundary), $|M|$ always denotes the $k$-dimensional Hausdorff measure of $M$. $W_{k-2}^{\bbS^{n}}(\wh{\del\Si})$ denotes the standard $(k-2)$-th quermassintegral of the closed convex hypersurface $\del\Si\sub\bbS^{n}$. We refer to \cref{QM} for a description via curvature integrals and for more information. Furthermore we define the {\it{spherical caps of radius $R$ around $e\in \bbS^{n}$}} by
\eq{C_{R}(e)=\{x\in \bar{\mathbb{B}}^{n+1}\cn |x-(R^{2}+1)^{\fr 12}e|= R\},\q R<\8,}
where we will drop the argument $e$, in cases where it is not relevant. Let $H$ denote the mean curvature of a hypersurface, i.e. $H=\s_{1}$.

The main results we will present are geometric inequalities relating some  particular quantities from \eqref{Geom-Quant} of convex free boundary hypersurfaces in the unit ball:
 
\begin{thm}\label{Main}
Let $n\geq 2$ and $\Si\sub \bar{\B}^{n+1}$ be a convex free boundary hypersurface in the unit ball. Then
 \eq{W_{n+1}(\wh\Si)=\fr{\om_{n}}{2(n+1)}} and for $0\leq k\leq n-1$ there holds
\eq{W_{n}(\wh\Si)\geq (f_{n}\circ f_{k}^{-1})(W_{k}(\wh\Si)),}
where $f_{k}=f_{k}(r)$ is the strictly increasing real function
\eq{f_{k}(r)=W_{k}(\wh{C_{r}}).}
Equality holds if and only if $\Si$ is a spherical cap or a flat disk. 
\end{thm}

\begin{rem}
The monotonicity of the functions $f_{k}$ is proven in \cref{Main-Proof}.
\end{rem}

We obtain the two following important special cases immediately:

\begin{cor}\label{Main-Cor}
Under the assumptions of \cref{Main}, for $n=2$ there holds the Minkowski type inequality
\eq{\fr{1}{6}\br{\int_{\Si}H+|\wh{\del\Si}|}\geq f_{2}\circ f_{1}^{-1}\br{\fr{1}{3}|\Si|}} 
and for $n=3$ we obtain a Willmore type inequality
\eq{\label{Willmore}\fr{1}{12}\br{\fr 13\int_{\Si}H^{2}+|\del\Si|}\geq f_{3}\circ f_{1}^{-1}\br{\fr 14|\Si|}, }
where in each case equality holds precisely on spherical caps or flat disks.
\end{cor}

Note that in \cite{Volkmann:/2016} Volkmann obtained an inequality similar to \eqref{Willmore} (with equality on caps and disks) in case $n=2$ for a much broader class of hypersurfaces using methods from geometric measure theory. In arbitrary higher dimensions, a Willmore type estimate for the convex case was deduced in \cite{LambertScheuer:09/2017}, however with equality only on flat disks. Hence \eqref{Willmore} is an extension of this result in dimension $n=3$.

The method of proof follows the nowadays classical method for proving geometric inequalities by employing monotonicity properties along and convergence of a suitable curvature flow. In the case of closed hypersurfaces, such methods were used in \cite{GuanLi:08/2009} to deduce Alexandrov-Fenchel inequalities and in \cite{HuiskenIlmanen:/2001} to prove the Riemannian Penrose inequality.  In other ambient spaces there is a zoo of variants of such inequalities, e.g. \cite{AndrewsChenWei:05/2018,De-LimaGirao:04/2016,GeWangWu:04/2014,HuLi:07/2018,WangXia:07/2014,WeiXiong:/2015} in the hyperbolic space, \cite{GiraoPinheiro:12/2017,MakowskiScheuer:11/2016,WeiXiong:/2015} in the sphere and \cite{BrendleHungWang:01/2016,GeWangWuXia:03/2015} in more general ambient spaces. 

The most suitable flows for this strategy are flows that decrease or increase a certain quantity, while preserving another one. Together with the knowledge that the flow converges to a well understood limit object, this implies the inequality. 

In the paper \cite{GuanLi:/2015}, Guan and Li constructed the mean curvature type flow 
\eq{\label{VPL}\dot{x}=(n-\s_{1}\ip{x}{\nu})\nu}
for closed and starshaped hypersurfaces in $\R^{n+1}$.
The Minkowski identity
\eq{\int_{M} n=\int_{M}{\s_{1}\ip{x}{\nu}}}
shows that the enclosed volume is preserved, while the second Minkowski identity
\eq{\int_{M}n\s_{1}=\int_{M}\s_{2}\ip{x}{\nu}}
shows that the surface area is decreasing. Together with suitable a priori estimates this gave a new flow approach to the isoperimetric inequality for starshaped hypersurfaces. Also in \cite{GuanLi:/2015}, a similar flow was considered in any space forms. Further generalisations of such flows were treated in \cite{GuanLi:/2018,GuanLiWang:09/2016,ScheuerXia:08/2017}. A suitably modified version of \eqref{VPL} for free boundary hypersurfaces in the unit ball will be treated in \cite{WangXia}. Possible generalisations also contain inverse type flows
\eq{\dot{x}=\br{\fr{c_{n,k}}{\fr{\s_{k}}{\s_{k-1}}}-\ip{x}{\nu}\nu},}
which have similar monotonicity properties due to the higher order Minkowski identities, cf. \cite{BrendleGuanLi:/}.

In the recent preprint \cite[Prop.~5.1]{WangXia:08/2017}, GW and CX have proven new Minkowski identities for immersed hypersurfaces $\Si$ with free boundary in the unit ball:
\eq{\int_{\Si}\s_{k-1}\ip{x}{e}=\fr{k}{n+1-k}\int_{\Si}\s_{k}\ip{X_{e}}{\nu},\q 1\leq k\leq n.}
Here $e\in \bbS^{n}$ is fixed and $X_{e}$ is the conformal Killing field
\eq{X_{e}=\ip{x}{e}x-\fr 12(|x|^{2}+1)e.}
The advantage of these Minkowski identities is that they do not contain any boundary integrals.

Motivated by these identities as well as the original flow \eqref{VPL} due to Guan and Li, it is natural to consider the following  inverse type flow for strictly convex hypersurfaces with free boundary in the unit ball:

\eq{\dot{x}=\br{\fr{\ip{x}{e}}{n\fr{\s_{n}}{\s_{n-1}}}-\ip{X_{e}}{\nu}}\nu}
with a suitably chosen direction $e\in\bbS^{n}$. We will prove that this curvature flow (with $e$ to be chosen in dependence of the initial data) will drive the strictly convex initial data in an infinite amount of time to a spherical cap. On the other hand we will show that the quantities in \eqref{Geom-Quant} have certain monotonicity properties, which will allow us to conclude the proof of \cref{Main}. The restriction to $\s_{n}$ is due to the technical obstruction, that for the general curvature function $\s_{k}/\s_{k-1}$ we could not prove the preservation of convexity along this flow. However, we still conjecture that all the lower order Alexandrov-Fenchel inequalities are true for convex free boundary hypersurfaces.

Finally we mention some previous results on curvature flows with free boundary in the unit ball. The classical mean curvature flow was considered in \cite{Stahl:08/1996,Stahl:06/1996}. The inverse mean curvature flow was treated in \cite{LambertScheuer:04/2016}, where it was shown that strictly convex initial data are driven to a flat perpendicular disk. We will also use several estimates for strictly convex hypersurfaces with free boundary in the unit ball, which were proven in this paper.

Our exposition is organised as follows.
In \cref{QM} we collect our conventions on geometric objects and we recall the definition of the quermassintegrals. In \cref{Sec:Flow} we prove the relevant properties of our flow and in \cref{Main-Proof} we finish the proof of \cref{Main}.

\section{Quermassintegrals}\label{QM}
\subsection*{Notation and basic definitions}\label{Sec:Basics}
In this section we state the basic conventions used in this paper. Let $M$ be a smooth manifold and $g$ be a Riemannian metric on $M$ with Levi-Civita connection $\n.$
Our convention for the Riemannian curvature tensor $\Rm$ is
\eq{\Rm(X,Y)Z=\n_{X}\n_{Y}Z-\n_{Y}\n_{X}Z-\n_{[X,Y]}Z} 
and the purely covariant version is
\eq{\Rm(X,Y,Z,W)=g(\Rm(X,Y)Z,W).}

Let 
\eq{x\cn M\ra \R^{n+1}}
be the smooth embedding of an $n$-dimensional closed and connected manifold. The induced metric of $x(M)$ is given by the pullback of the ambient Euclidean metric $\bar g$ by $x$.
 The second fundamental form $h$ of the embedding $x$ is given by the Gaussian formula
\eq{\bar\n_{X}Y=\n_{X}Y- h(X,Y)\nu.}
The Weingarten operator is defined via
\eq{g(\cW(X),Y)=h(X,Y)}
and the Weingarten equation says that
\eq{\label{Weingarten}\bar{\n}_{X}\nu=\cW(X).}  
Finally, we have Gauss equation,
\eq{\label{GaussEq}\Rm(W,X,Y,Z)&=\br{h(W,Z)h(X,Y)-h(W,Y)h(X,Z)}.}
	
\begin{rem}\label{Rem-coordinates}
We will simplify the notation by using the following shortcuts occasionally:
\begin{enumerate}[(i)]

\item When dealing with complicated evolution equations of tensors, we will use a local frame to express tensors with the help of their components, i.e. for a tensor field $T\in \cT^{k,l}(M)$, the expression $T^{i_{1}\dots i_{k}}_{j_{1}\dots j_{l}}$ denotes  
\eq{T^{i_{1}\dots i_{k}}_{j_{1}\dots j_{l}}=T(e_{j_{1}},\dots,e_{j_{l}},\e^{i_{1}},\dots \e^{i_{k}}),}
where $(e_{i})$ is a local frame and $(\e^{i})$ its dual coframe. 
\item The coordinate expression for  the $m$-th covariant derivative of a $(k,l)$-tensor field $T$ is 
\eq{\n^m T=\br{T^{i_1\dots i_k}_{j_1\dots j_l;j_{l+1}\dots j_{l+m}}},}
where indices appearing after the semi-colon denote the covariant derivatives.
\item For convenience the components of the Weingarten map $\cW$ are denoted by $(h^{i}_{j})$.
\end{enumerate}
\end{rem}

\subsection*{Quermassintegrals in {\texorpdfstring{$\bbR^{n+1}$}{}}} In order to compare our new geometric integrals for free boundary  hypersurfaces
with the quermassintegrals for closed  hypersurfaces in a space form, we start to introduce the quermassintegrals  in  $\bbR^{n+1}$. Let $\Omega$ be a convex body 
(non-empty, compact, convex set) in $\bbR^{n+1}$. The  quermassintegrals $W_{k}^{\bbR^{n+1}}$ are defined by

\eq{
W_0^{\bbR^{n+1}}(\Omega)&= |\Omega|, \quad W_1^{\bbR^{n+1}}(\Omega)=\frac{1}{n+1}|\partial \Omega|,\\
W_{k}^{\bbR^{n+1}}(\Omega)&= \frac{1}{n+1}\binom{n}{k-1}^{-1}\int_{\partial \Omega}\s_{k-1} dA\\
&=  \frac{1}{n+1}\int_{\partial \Omega} H_{k-1} dA, \quad
\quad  2\le k\le n+1.
}
In particular, 
\eq{W_{n+1}^{\bbR^{n+1}}(\Omega)=\frac{\omega_{n}}{n+1}.}
The first variation formula of $W_{k}^{\bbR^{n+1}}$ is as follows.  Consider a family of bounded convex bodies $\{\Omega_t\}$ in $\bbR^{n+1}$  whose boundary  
$\partial \Omega_t$ evolving by a normal variation with speed function $f$, then
\begin{eqnarray}\label{1st-var_a}
\frac{d}{dt}W_{k}^{\bbR^{n+1}}(\Omega_t)=\frac {n+1-k} {n+1} \int_{\del K_t} H_k f dA_t, \q 0\leq k\leq n+1.
\end{eqnarray}

\subsection*{Quermassintegrals in {\texorpdfstring{$\bbS^n$}{}}}
For a convex body (non-empty, compact, convex set) $K\subset \bbS^n$ with smooth boundary, the quermassintegrals $W_{k}^{\bbS^n}$ are inductively defined by 

\eq{
W_0^{\bbS^n}(K)&= |K|, \quad W_1^{\bbS^n}(K)=\frac{1}{n}|\partial K|,\\
W_{k}^{\bbS^n}(K)&= \frac{1}{n}\binom{n-1}{k-1}^{-1}\int_{\partial K}\s^{\bbS^n}_{k-1} dA +\frac{k-1}{n-k+2}W_{k-2}^{\bbS^n}(K)\\
&= \frac{1}{n}\int_{\partial K} H^{\bbS^n}_{k-1} dA +\frac{k-1}{n-k+2}W_{k-2}^{\bbS^n}(K), \quad 
\quad  2\le k\le n,
}
where by $\s_{k}^{\bbS^{n}}$ we denote the $k$-th elementary symmetric polynomial, evaluated at the $n-1$ principal curvatures of the hypersurface $\del K\sub\bbS^{n}$ and 
\eq{H_{k}^{\bbS^{n}}=\fr{1}{\binom{n-1}{k}}\s_{k}^{\bbS^{n}},\q 0\leq k\leq n-1,\q H_{n}^{\bbS^{n}}=\s_{n}^{\bbS^{n}}=0.}
In particular, 
\eq{W_n^{\bbS^n}(K)=\frac{\omega_{n-1}}{n}} due to the spherical Gauss-Bonnet-Chern's Theorem, cf. \cite{Solanes:/2006}.

The first variation formula of $W_{k}^{\bbS^n}$ is as follows.  Consider a family of bounded convex bodies $\{K_t\}$ in $\bbS^n$  whose boundary  
$\partial K_t$ evolves by a normal variation with speed function $f$, we have similarly 
\eq{\label{1st-var}
\frac{d}{dt}W_{k}^{\bbS^n}(K_t)=\frac{1}{\binom{n}{k}}\int_{\del K_t} \sigma_k^{\bbS^n} f dA_t=\frac {n-k} n \int_{\del K_t} H_k^{\bbS^n} f dA_t,\q 0\leq k\leq n.
}
This can be proved by a simple induction argument; see for example \cite[Prop.~3.1]{WangXia:07/2014} for a similar deduction in the hyperbolic case.

\subsection*{Related quantities for free boundary hypersurfaces in {\texorpdfstring{$\bar\B^{n+1}$}{}}}
Let $\Si\subset \bar \B^{n+1}$ be a smooth convex, embedded hypersurface with free boundary $\partial \Si\subset \bbS^n$.
Denote by $\wh{\partial \Si}$ the convex body in $\bbS^n$ enclosed by $\partial \Si$. Denote by $\wh{\Si}$ the enclosed convex domain in $\bar \B^{n+1}$
which contains $\wh{\partial \Si}$.

We define $W_k$ for $\wh{\Si}\subset \bar \B^{n+1}$ as follows.
\eq{\label{Def-QM}W_{0}(\wh\Si)&=|\wh\Si|,\quad W_{1}(\wh\Si)= \frac 1{n+1} |\Si|,\\
	W_{k}(\wh\Si)&=   \frac 1{n+1} 
	\binom{n}{k-1}^{-1}
	\bigg\{\int_{\Si}\s_{k-1}dA+\binom{n}{k-2}W_{k-2}^{\bbS^{n}}(\wh{\del\Si})\bigg\}\\
	&= \frac 1 {n+1} \int_\Sigma H_{k-1} dA + \frac 1 {n+1} \frac {k-1}{n-k+2}W_{k-2}^{\bbS^{n}}(\wh{\del\Si}),
	\quad  2\leq k\leq n+1.}

It is the following variational formula which motivates our definition of the $W_{k}$.

\begin{prop}\label{Lemma-Wk}
Let $\{\Si_t\}$ be a family of smooth, embedded hypersurfaces with free boundary in $\bar \B^{n+1}$, given by the embeddings
$x(\cdot, t): \bar\B^n\to \bar\B^{n+1}$, which evolve by 
\begin{eqnarray*}
\del_t x(\cdot, t)= f(\cdot, t)\nu(\cdot, t).
\end{eqnarray*}
Then there holds
\eq{\label{Wk}\del_{t}W_{k}(\wh{\Si_t})=\frac {n+1-k}{n+1} \int_{\Si_{t}}H_{k}f~dA_{t}, \quad  \q \fa 0\leq k\leq n+1.
} In particular,
\eq{\label{Wk-top}\del_{t}W_{n+1}(\wh{\Si_t})=0.
}
\end{prop}

\pf{The following general evolution equations are valid
\eq{\dot{g}_{ij}=2fh_{ij},} 
\eq{\dot{h}^{i}_{j}=-{f_{;}}_{j}^{i}-fh^{i}_{m}h^{m}_{j}.}
 The variational formulas for $W_0$ and $W_1$ are well-known. 
Consider the case $k\ge 2$. We prove the following equivalent formula of \eqref{Wk}.
\begin{equation}
 \label{Wk-a}\partial_t\bigg\{\int_{\Si}\s_{k-1}dA+\binom{n}{k-2}W_{k-2}^{\bbS^{n}}(\wh{\del\Si})
 \bigg\} =\int_{\Si_{t}}k\s_{k}f~dA_{t}.
\end{equation}
We calculate, using
\eq{&\fr{\del\s_{k-1}}{\del h^{j}_{i}} h^{i}_{m}h^{m}_{j}=\s_{1}\s_{k-1}-k\s_{k},\quad 2\le k\le n,\\
&\fr{\del\s_{n}}{\del h^{j}_{i}}h^{i}_{m}h^{m}_{j}=\s_{1}\s_{n},}
that
\eq{\del_{t}\int_{\Si_{t}}\s_{k-1}~dA_{t}&=\int_{\Si_{t}}\s_{k-1}f\s_{1}~dA_{t}-\int_{\Si_{t}}\fr{\del \s_{k-1}}{\del h^{j}_{i}}{f_{;}}^{i}_{j}~dA_{t}-\int_{\Si_{t}}\fr{\del\s_{k-1}}{\del h^{j}_{j}}h^{i}_{m}h^{m}_{j}f~dA_{t}\\
					&=\left\{\begin{array}{lll}\int_{\Si_{t}}k\s_{k}f~dA_{t}-\int_{\del \Si_{t}}\fr{\del \s_{k-1}}{\del h^{j}_{i}}{f_{;}}_{j}\mu^{i}~dA_{t}, \quad &2\le k\le n\\
					-\int_{\del \Si_{t}}\fr{\del \s_{k-1}}{\del h^{j}_{i}}{f_{;}}_{j}\mu^{i}~dA_{t},\quad &k=n+1,
					\end{array}\right.}
where we used that the $\s_{k}$ are divergence-free in $\R^{n+1}$.
Let $\{e_{I}\}_{2\leq I\leq n}$ be an orthonormal frame for $\del\Si$, so that $\{e_{1}=\mu, (e_{I})_{2\leq I\leq n}\}$ is an orthonormal frame for $\Si$ along $\del\Si$.
Since $\mu$ is a principal direction due to the free boundary condition, there holds \eq{\fr{\del \s_{k-1}}{\del h^{J}_{\mu}}=0,\quad 2\leq J\leq n.}
Differentiation of $\ip{\nu}{\mu}=0$ with respect to time gives
\eq{\left<-\nabla f, \mu\right>+f h_{\nu\nu}^{\bbS^n}=0,} which implies
\eq{\label{BD-f}\n_{\mu}f=f.}
It follows that
\eq{\int_{\del\Si_{t}}\fr{\del \s_{k-1}}{\del h^{j}_{i}}{f_{;}}_{j}\mu^{i}~dA_{t}=\int_{\del \Si_{t}}\fr{\del \s_{k-1}}{\del h^{\mu}_{\mu}}f~dA_{t}.}

On the other hand, the flow of $x(\cdot, t)$ induces a normal flow in $\bbS^n$ with speed $f$. It follows from \eqref{1st-var} that
\eq{\del_{t}W_{k-2}^{\bbS^{n}}(\wh{\del \Si_{t}})=\fr{1}{\binom{n}{k-2}}\int_{\del \Si_{t}}\s_{k-2}^{\bbS^{n}}f~dA_{t}.}
Combining these equalities, also having in mind that
\eq{\fr{\del\s_{k-1}}{\del h^{\mu}_{\mu}}=\s_{k-2}^{\bbS^{n}},} gives \eqref{Wk-a}, and hence \eqref{Wk}.}

\begin{rem} The variation formulas of  geometric integrals $W_k$, which we define for free boundary hypersurfaces,
\eq{
 \del_{t}W_{k}(\wh{\Si_t})=\frac{n+1-k}{n+1}\int_{\Si_{t}}H_{k}f~dA_{t} \q \fa 1\leq k\leq n, \label{eq_a}
 }
 are similar to those for quermassintegrals for closed hypersurfaces in a space form. Therefore we believe that these geometric integrals
 are correct counterparts to
the classical quermassintegrals. In fact, for $W_{n+1}$ we also have a very similar form, see the next Proposition. 
It is an interesting question whether one can establish an integral geometric theory for 
free boundary hypersurfaces in a ball.
\end{rem}

\begin{prop}\label{GB}
Let $n\geq 2$ and $\Si\sub \bar{\B}^{n+1}$ be a convex free boundary hypersurface in the unit ball. Then
\eq{W_{n+1}(\wh{\Si})=\fr{1}{2} \frac{\omega_n}{n+1}.
}
\end{prop}

The proof will be given in Appendix A. In fact, we will give a precise formula for $W_{n+1}$ for any free boundary hypersurface
in terms of the Euler characteristic.

\begin{rem}
For $n=1$ there holds
\eq{W_2(\wh\Si)=\fr{1}{2}\br{\int_{\Si} \kappa ds +|\wh{\partial \Si}|}=\fr{\pi}{2}.}

For $n=2$ we have
\eq{W_2(\wh\Si)=\fr 16\br{\int_{\Si} H dA + |\wh{\partial \Si}|}} and \eq{W_3(\wh\Si)=\fr 13\br{\int_{\Si} \s_2 dA +|\partial \Si|}= \fr{2\pi}{3},}


\end{rem}


\section{The curvature flow}\label{Sec:Flow}
Throughout this section $\Si$ is as in \cref{Main} and in addition strictly convex. 
In order to prove the proposed geometric inequality, we use the following curvature flow with free perpendicularity condition on $\bbS^{n}$,
\eq{\label{Flow}\dot{x}=f\nu&\equiv \br{\fr{\ip{x}{e}}{F}-\ip{X_{e}}{\nu}}\nu\\
		x(0,\bar{\B}^{n})&=\Si_{0}:=\Si,}
where $e\in \wh{\del\Si}$ is chosen such that $\wh{\del\Si}$ lies in the open hemisphere
\eq{\cH=\{x\in \bbS^{n}\cn \ip{x}{e}>0\}} 
and where $X_{e}$ is the conformal Killing field
\eq{X_{e}=X_{e}(x)=\ip{x}{e}x-\fr 12(|x|^{2}+1)e.}
In particular we will use the curvature function
\eq{\label{F}F=n\fr{\s_{n}}{\s_{n-1}}}
and denote the flow hypersurfaces by $\Si_{t}=x(t,\bar{\B}^{n})$.
The normalisation is chosen such that
\eq{F(1,\dots,1)=1.}

\subsection*{Monotonicity}

The main motivation to study the flow \eqref{Flow} comes from its monotonicity properties.
From \cref{Lemma-Wk} we obtain:

\begin{prop}\label{Monotone}
Along the flow \eqref{Flow}, the quantity $W_{n}(\wh{\Si_{t}})$ is preserved and $W_{k}(\wh{\Si_{t}})$ is increasing for all $0\leq k\leq n-1$. In case $0\leq k\leq n-1$, $W_{k}(\wh{\Si_{t}})$ is constant if and only if $\Si_{t}$ is a spherical cap.
\end{prop}

\pf{
Expressing \eqref{Wk-a} for $1\leq k\leq n$ explicitly gives
\eq{(n+1)\binom{n}{k-1}\del_{t}W_{k}(\wh{\Si_{t}})&=\fr{k}{n}\int_{\Si_{t}}\br{\fr{\s_{k}\s_{n-1}}{\s_{n}}\ip{x}{e}-n\s_{k}\ip{X_{e}}{\nu}}\\
			&\geq k\int_{\Si_{t}}\br{\fr{n-k+1}{k}\s_{k-1}\ip{x}{e}-\s_{k}\ip{X_{e}}{\nu}}\\
			&=0,}
due to Newton's inequality and \cite[Prop.~5.1]{WangXia:08/2017}.
If $k=n$, we have equality. For $k=0$ we have
\eq{\del_{t}W_{0}(\wh{\Si_{t}})&=\int_{\Si_{t}}\br{\fr{\ip{x}{e}}{\fr{n\s_{n}}{\s_{n-1}}}-\ip{X_{e}}{\nu}}\geq \int_{\Si_{t}}\br{\fr{n\ip{x}{e}}{H}-\ip{X_{e}}{\nu}}\geq 0, }
due to the Heintze-Karcher-type inequality in \cite[Thm.~5.2]{WangXia:08/2017}. This proves the first claim, while the second claim is due to the equality case of Newton's inequality.
}

In the next part of this section we prove the smooth convergence of the solution of \eqref{Flow} to a spherical cap. The previous monotonicity properties will then finish the proof of \cref{Main}.

\subsection*{Barriers}
Contrary to a purely expanding inverse curvature flow, as it was considered in \cite{LambertScheuer:04/2016} and applied in \cite{LambertScheuer:09/2017} to give an estimate for a Willmore type quantity for convex free boundary hypersurfaces, our flow \eqref{Flow} will stay away from the minimal disk and hence, at least in principle, natural singularities are avoided. This nice feature is due to the constraining terms $\ip{x}{e}$ and $\ip{X_{e}}{\nu}$, which force the flow to remain between two spherical barriers.

The short-time existence is established by standard PDE theory, since due to our choice of $e$ there holds $\ip{x}{e}_{|\Si_{0}}>0$ and hence the operator is uniformly parabolic, compare \cite{Marquardt:/2012} for instance. Let $T^{*}$ denote the maximal time of smooth existence of a solution to \eqref{Flow}. This implies that all flow hypersurfaces are strictly convex up to $T^{*}$ due to the positivity of $F$ up to $T^{*}$. 
The family of spherical caps lying entirely in the closed half-ball $\bar{\mathbb{B}}^{n+1}_{+}$ is given by
\eq{C_{R}(e)=\{x\in \bar{\mathbb{B}}^{n+1}\cn |x-(R^{2}+1)^{\fr 12}e|= R\},\q R<\8.}
These form the natural barriers of \eqref{Flow}:

\begin{lemma}\label{Barriers}
Suppose for two radii $R_{1}<R_{2}$, $M_{0}$ satisfies
\eq{\Si_{0}\sub \wh{C_{R_{2}}(e)}\bs \wh{C_{R_{1}}(e)},     }
then this is preserved along the flow. In particular the height function $\ip{x}{e}$ satisfies a priori bounds
\eq{\e\leq \ip{x}{e}\leq 1-\e}
 and the normal vector points uniformly downwards,
 \eq{\ip{\nu(t,\xi)}{e}\leq -\ti\e\q\fa (t,\xi)\in [0,T^{*})\x \bar{\B}^{n},}
for some suitable $\e,\ti\e\in(0,1)$, which only depend on the initial datum. 
\end{lemma}

\pf{
An elementary calculation shows that the spherical caps $C_{R}(e)$ are static solutions of \eqref{Flow}, i.e. $f=0$ along $C_{R}(e)$. A simple comparison principle shows that the caps are barriers, since the avoidance principle holds due to the free boundary condition. The height estimate follows immediately. The final claim is the statement of \cite[Lemma~11]{LambertScheuer:04/2016}.
}

\subsection*{Evolution equations}
We collect the relevant evolution equations, which will allow us to handle the curvature flow. In order to avoid confusion with tensor indices, here we abbreviate
\eq{X=X_{e},}
keeping in mind that $X$ depends on $e$. We need the following lemma.

\begin{lemma}
There hold
\eq{\ip{X}{\nu}_{;i}=\ip{x_{;i}}{e}\ip{x}{\nu}-\ip{x_{;i}}{x}\ip{e}{\nu}+h^{k}_{i}\ip{X}{x_{;k}}}
and
\eq{\ip{X}{\nu}_{;ij}&=h^{k}_{i;j}\ip{X}{x_{;k}}+\ip{x}{e}h_{ij}-h^{k}_{i}h_{kj}\ip{X}{\nu}-g_{ij}\ip{e}{\nu}\\
				&\hp{=}+h^{k}_{j}\br{\ip{x_{;i}}{e}\ip{x}{x_{;k}}-\ip{x_{;i}}{x}\ip{e}{x_{;k}}}\\
				&\hp{=}+h^{k}_{i}\br{\ip{x_{;j}}{e}\ip{x}{x_{;k}}-\ip{x_{;j}}{x}\ip{e}{x_{;k}}}}
\end{lemma}

\pf{
\eq{\ip{X}{\nu}_{;i}&=\ip{X_{;i}}{\nu}+\ip{X}{\nu_{;i}}\\
			&=\ip{x_{;i}}{e}\ip{x}{\nu}-\ip{x_{;i}}{x}\ip{e}{\nu}+h^{k}_{i}\ip{X}{x_{;k}}.}
\eq{\ip{X}{\nu}_{;ij}&=-h_{ij}\ip{e}{\nu}\ip{x}{\nu}+\ip{x_{;i}}{e}h^{k}_{j}\ip{x}{x_{;k}}\\
			&\hp{=}+h_{ij}\ip{\nu}{x}\ip{e}{\nu}-g_{ij}\ip{e}{\nu}-\ip{x_{;i}}{x}h^{k}_{j}\ip{e}{x_{;k}}\\
			&\hp{=}+h^{k}_{i;j}\ip{X}{x_{;k}}+h^{k}_{i}\ip{\ip{x_{;j}}{e}x+\ip{x}{e}x_{;j}-\ip{x_{;j}}{x}e}{x_{;k}}-h^{k}_{i}h_{kj}\ip{X}{\nu}\\
			&=h^{k}_{j}\br{\ip{x_{;i}}{e}\ip{x}{x_{;k}}-\ip{x_{;i}}{x}\ip{e}{x_{;k}}}-g_{ij}\ip{e}{\nu}+h^{k}_{i;j}\ip{X}{x_{;k}}\\
			&\hp{=}+h^{k}_{i}\br{\ip{x_{;j}}{e}\ip{x}{x_{;k}}-\ip{x_{;j}}{x}\ip{e}{x_{;k}}}+h_{ij}\ip{x}{e}-h_{i}^{k}h_{kj}\ip{X}{\nu}.}
}

Let 
\eq{\cL=\del_{t}-\fr{\ip{x}{e}}{F^{2}}F^{ij}\n^{2}_{ij}-\ip{X}{\n},}
where $F$ may either be understood to depend on the Weingarten operator, $F=F(h^{i}_{j})$ or on the second fundamental form and the metric, $F=F(h_{ij},g_{ij}).$
There holds
\eq{F^{i}_{j}\equiv\fr{\del F}{\del h^{j}_{i}}=g_{kj}F^{ki}\equiv g_{kj}\fr{\del F}{\del h_{ki}},}
cf. \cite[Ch.~2]{Gerhardt:/2006}.
The following specific evolution equations are valid. 

\begin{prop}\label{Ev-Eq}
For general 1-homogeneous $F$, along \eqref{Flow} there hold
\enum{
\item 
	\eq{\label{Ev-Height}\cL\ip{x}{e}=\fr{2}{F}\ip{x}{e}\ip{\nu}{e}-\ip{X}{e}.}
 \eq{\label{BD-Height}\n_{\mu}\ip{x}{e}=\ip{x}{e}\q\mbox{along}~\del\Si_{t}.} 
\item  \eq{\label{Ev-F}\cL F&=-\fr{2\ip{x}{e}}{F^{3}}F^{ij}F_{;i}F_{;j}+\fr{2}{F^{2}}F^{ij}\ip{x_{;i}}{e}F_{;j}\\
	&\hp{=}+\br{1-F^{ij}g_{ij}}\ip{e}{\nu}+F\ip{x}{e}\br{1-\fr{F^{ij}h^{k}_{i}h_{kj}}{F^{2}}}.}
	\eq{\label{BD-F}\n_{\mu}F=0\q\mbox{along}~\del\Si_{t}.}

\item 

\eq{\label{Ev-H}\cL{H}&=\fr{\ip{x}{e}}{F^{2}}g^{ij}F^{kl,rs}h_{kl;i}h_{rs;j}+\br{\fr{F^{kl}h_{km}h^{m}_{l}}{F^{2}}+1}\ip{x}{e}H-\fr{2\ip{x}{e}}{F}\|\cW\|^{2}\\
		&\hp{=}+\br{\fr{H}{F}-n}\ip{\nu}{e}+\fr{2}{F^{2}}g^{ij}\ip{x_{;i}}{e}F_{;j}-\fr{2\ip{x}{e}}{F^{3}}\|\n F\|^{2}.}
\eq{\label{BD-H}\n_{\mu}H\leq 0\q\mbox{along}~\del\Si_{t},}
provided $F$ is concave.
}

\end{prop}

\pf{
(i) (a)
\eq{\del_{t}\ip{x}{e}-\fr{\ip{x}{e}}{F^{2}}F^{ij}\ip{x}{e}_{;ij}&=\fr{\ip{x}{e}\ip{\nu}{e}}{F}-\ip{X}{\nu}\ip{\nu}{e}+\fr{\ip{x}{e}\ip{\nu}{e}}{F}\\
					&=\fr{2}{F}\ip{x}{e}\ip{\nu}{e}-\ip{X}{\nu}\ip{\nu}{e}\\
					&=\fr{2}{F}\ip{x}{e}\ip{\nu}{e}-\ip{X}{e}+{\ip{x}{e}_{;}}^{i}\ip{X}{x_{;i}}.}
					
(b)~was deduced in \cite[Lemma~5]{LambertScheuer:04/2016}.

(ii) (a) \eq{\dot{F}&=F^{j}_{i}\dot{h}^{i}_{j}\\
		&=-F^{j}_{i}{\br{\fr{\ip{x}{e}}{F}-\ip{X}{\nu}}_{;}}_{j}^{i}-F^{j}_{i}\br{\fr{\ip{x}{e}}{F}-\ip{X}{\nu}}h^{i}_{k}h^{k}_{j}\\
		&=-F^{ij}\br{\fr{-h_{ij}\ip{\nu}{e}}{F}-\fr{\ip{x}{e}}{F^{2}}F_{;ij}+\fr{2\ip{x}{e}}{F^{3}}F_{;i}F_{;j}-\fr{2\ip{x_{;i}}{e}}{F^{2}}F_{;j}}\\
		&\hp{=}+{F_{;}}^{k}\ip{X}{x_{;k}}+\ip{x}{e}F-F^{ij}h^{k}_{i}h_{kj}\ip{X}{\nu}-F^{ij}g_{ij}\ip{e}{\nu}\\
		&\hp{=}-\fr{\ip{x}{e}}{F}F^{ij}h^{k}_{i}h_{kj}+F^{ij}h_{i}^{k}h_{kj}\ip{X}{\nu}\\
		&=\fr{\ip{x}{e}}{F^{2}}F^{ij}F_{;ij}-\fr{2\ip{x}{e}}{F^{3}}F^{ij}F_{;i}F_{;j}+\fr{2}{F^{2}}F^{ij}\ip{x_{;i}}{e}F_{;j}+{F_{;}}^{k}\ip{X}{x_{;k}}\\
		&\hp{=}+\br{1-F^{ij}g_{ij}}\ip{e}{\nu}+F\ip{x}{e}\br{1-\fr{F^{ij}h^{k}_{i}h_{kj}}{F^{2}}}.} 
		
(b)~There holds, using an orthonormal frame $(e_{I})_{2\leq I\leq n}$ for $T_{x}(\del\Si_{t})$,
\eq{\n_{\mu}\ip{X}{\nu}_{|\del \Si_{t}}&=\n_{\mu}\br{\ip{x}{e}\ip{x}{\nu}-\fr 12(|x|^{2}+1)\ip{e}{\nu}}_{|\del \Si_{t}}\\
						&=\ip{\mu}{e}\n_{\mu}\ip{x}{\nu}-\ip{e}{\nu}-\n_{\mu}\ip{e}{\nu}\\
						&=-\ip{e}{\nu}+\ip{\ip{\mu}{e}\mu-e}{\cW(\mu)}\\
						&=-\ip{e}{\nu}+\ip{\ip{\mu}{e}\mu-e}{\ip{\cW(\mu)}{\mu}\mu+\sum_{I=2}^{n}\ip{\cW(\mu)}{e_{I}}e_{I}}\\
						&=-\ip{e}{\nu}\\
						&=\ip{X}{\nu}_{|\del \Si_{t}}.}
From \eqref{BD-Height} and \eqref{BD-f} we obtain 
\eq{\n_{\mu}F=0.} 
(iii) (a)
\eq{\label{Ev-H-1}\dot{H}&=-g^{ij}f_{;ij}-f\|\cW\|^{2}\\
		&=g^{ij}\br{\ip{X}{\nu}-\fr{\ip{x}{e}}{F}}_{;ij}-f\|\cW\|^{2}\\
		&={H_{;}}^{k}\ip{X}{x_{;k}}+H\ip{x}{e}-\|\cW\|^{2}\ip{X}{\nu}-n\ip{e}{\nu}+\fr{H}{F}\ip{\nu}{e}\\
		&\hp{=}+\fr{2}{F^{2}}g^{ij}\ip{x_{;i}}{e}F_{;j}+\fr{\ip{x}{e}}{F^{2}}g^{ij}F_{;ij}-\fr{2\ip{x}{e}}{F^{3}}\|\n F\|^{2}-f\|\cW\|^{2}\\
		&={H_{;}}^{k}\ip{X}{x_{;k}}+H\ip{x}{e}-\fr{\ip{x}{e}}{F}\|\cW\|^{2}-n\ip{e}{\nu}+\fr{H}{F}\ip{\nu}{e}\\
		&\hp{=}+\fr{2}{F^{2}}g^{ij}\ip{x_{;i}}{e}F_{;j}+\fr{\ip{x}{e}}{F^{2}}g^{ij}F^{kl}h_{kl;ij}+\fr{\ip{x}{e}}{F^{2}}g^{ij}F^{kl,rs}h_{kl;i}h_{rs;j}\\
		&\hp{=}-\fr{2\ip{x}{e}}{F^{3}}\|\n F\|^{2}.}
Now we have to interchange the pairs $(h,k)$ with $(i,j)$ in $h_{kl;ij}$. Due to the Codazzi equations, the Ricci identities and the Gauss equation there holds
\eq{
h_{kl;ij}&=h_{ki;lj}\\
                &=h_{ki;jl}+{R_{ljk}}^ah_{ai}+{R_{lji}}^a h_{ka}\\
                &=h_{ij;kl}+{R_{ljk}}^ah_{ai}+{R_{lji}}^a h_{ka}\\
                         &=h_{ij;kl}+(h_{la}h_{jk}-h_{lk}h_{ja})h_{i}^a+(h_{la}h_{ji}-h_{li}h_{ja})h_{k}^a
             }
and thus
\eq{g^{ij}F^{kl}h_{kl;ij}=F^{kl}H_{;kl}+F^{kl}h_{la}h^{a}_{k}H-F\|\cW\|^{2}.}
Inserting into \eqref{Ev-H-1} gives the proposed evolution equation.

(b)~For an orthonormal frame $(e_{I})_{2\leq I\leq n}$ of $T_{x}(\del\Si_{t})$, such that the second fundamental form of $\Si$ is diagonal with respect to $(\nu,e_{I})_{2\leq I\leq n}$, along $\del\Si_{t}$ we calculated in \cite[Lemma~2]{LambertScheuer:04/2016}, that
\eq{ h_{IJ;\mu}=-h_{IJ}+h_{\mu\mu}\ip{e_{I}}{e_{J}}.}
From \eqref{BD-F} we obtain
\eq{0=\n_{\mu}F=\fr{\del F}{\del h}(\n_{\mu}h)=\sum_{I=2}^{n}F^{II}h_{II;\mu}+F^{\mu\mu}h_{\mu\mu;\mu}}
and thus
\eq{\n_{\mu}H&=h_{\mu\mu;\mu}+\sum_{I=2}^{n}h_{II;\mu}\\
			&=-\sum_{I=2}^{n}\fr{F^{II}}{F^{\mu\mu}}h_{II;\mu}-\sum_{I=2}^{n}h_{II}+(n-1)h_{\mu\mu}\\
			&=\sum_{I=2}^{n}\fr{F^{II}}{F^{\mu\mu}}h_{II}-\sum_{I=2}^{n}\fr{F^{II}}{F^{\mu\mu}}h_{\mu\mu}-\sum_{I=2}^{n}h_{II}+(n-1)h_{\mu\mu}\\
			&=\sum_{I=2}^{n}\fr{1}{F^{\mu\mu}}\br{F^{II}-F^{\mu\mu}}\br{h_{II}-h_{\mu\mu}},}
which is non-positive due to the concavity of $F$.
		}

\subsection*{Curvature estimates}
We continue the a priori estimates for \eqref{Flow}. 

\begin{prop}
Let $F$ be given by \eqref{F}. Then along \eqref{Flow} there hold: 
\enum{
\item The function
\eq{\p(t)=\min_{\bar{\B}^{n}}F(t,\cdot)}
is non-decreasing and hence
\eq{F(t,\xi)\geq \min_{\bar{\B}^{n}}F(0,\cdot)\q\fa (t,\xi)\in [0,T^{*})\x\bar{\B}^{n}.}
\item There exists a constant $c>0$, depending only on the initial data, such that 
\eq{F\leq c.}
}
\end{prop}

\pf{
(i)~
Considering $F$, and also $\s_{k}$, as a function of the principal curvatures $\ka_{i}$ and denoting
\eq{F^{i}=\fr{\del F}{\del\ka_{i}},}
we compute
\eq{F^{i}=n\fr{\s_{n}^{i}}{\s_{n-1}}-n\fr{\s_{n}\s_{n-1}^{i}}{\s_{n-1}^{2}}.}
Due to the relation
\eq{\s_{k}=\s_{k}^{i}\ka_{i}+\s_{k+1}^{i}\q\fa 0\leq k\leq n}
 we obtain
\eq{\sum_{i=1}^{n}\s_{k}^{i}=\sum_{i=1}^{n}\br{\s_{k-1}-\s_{k-1}^{i}\ka_{i}}=(n-k+1)\s_{k-1}.}
Hence
\eq{\sum_{i=1}^{n} F^{i}=n-2n\fr{\s_{n}\s_{n-2}}{\s_{n-1}^{2}}\geq 1}
due to the Newton-Maclaurin inequalities.
Furthermore there holds for $1\leq k\leq n$:
\eq{\s_{1}\s_{k}-(k+1)\s_{k+1}&=\sum_{i=1}^{n}\br{\ka_{i}\s_{k}-\ka_{i}\s_{k+1}^{i}}=\sum_{i=1}^{n}\ka_{i}^{2}\s_{k}^{i}}
and thus we compute
\eq{\sum_{i=1}^{n}F^{i}\ka_{i}^{2}=\fr{n}{\s_{n-1}}\s_{1}\s_{n}-n\fr{\s_{n}}{\s_{n-1}^{2}}(\s_{1}\s_{n-1}-n\s_{n})=n^{2}\fr{\s_{n}^{2}}{\s_{n-1}^{2}}=F^{2}.}
At the boundary there holds
\eq{\n_{\mu}F=0}
and hence at minimising points for $F$, the right hand side of its evolution \eqref{Ev-F} is non-negative and thus $\p$ is non-decreasing, taking \cref{Barriers} into account.

(ii)~To estimate $F$ from above we define
\eq{\p=\log F-\al\ip{x}{e}}
for some positive $\al$ which will be determined later.
We calculate the evolution equation
\eq{\cL\p&=\fr{1}{F}\cL F+\fr{\ip{x}{e}}{F^{2}}F^{ij}(\log F)_{;i}(\log F_{;j})-\al\cL\ip{x}{e}\\
		&=-\fr{\ip{x}{e}}{F^{2}}F^{ij}(\log F)_{;i}(\log F)_{;j}+\fr{2}{F^{2}}F^{ij}\ip{x}{e}_{;i}(\log F)_{;j}\\
	&\hp{=}+\fr{1}{F}\br{1-F^{ij}g_{ij}}\ip{e}{\nu}-\al\br{\fr{2}{F}\ip{x}{e}\ip{\nu}{e}-\ip{X}{e}}.}
On the boundary there holds
\eq{\n_{\mu}\p=-\al\ip{x}{e}<0.}
Hence maximal values of $\p$ are attained in the interior and at such we have
\eq{0&\leq \cL\p\leq \fr{1}{F^{2}}\br{-\al^{2}\ip{x}{e}+2\al}F^{ij}\ip{x}{e}_{;i}\ip{x}{e}_{;j}+\fr{c+2\al}{F}+\al\ip{X}{e}.}
The quantity
\eq{\ip{X}{e}&=\ip{x}{e}^{2}-\fr 12(|x|^{2}+1)=-\fr 12\br{|x|^{2}-\ip{x}{e}^{2}}-\fr 12\br{1-\ip{x}{e}^{2}}\leq -c_{\e}, }
due to \cref{Barriers}. Hence, picking $\al$ large enough in dependence of the previous a priori bounds, we obtain a contradiction if $F$ is too large.
}

To finish the curvature estimates, we use the evolution of the mean curvature $H=\s_{1}$ to prove that it is bounded. This is sufficient due to the convexity.

\begin{prop}
If $F$ is given by $\eqref{F}$, then the mean curvature, and hence all principal curvatures, are a priori bounded by the initial data.
\end{prop}

\pf{
In the evolution equation \eqref{Ev-H} we use the Cauchy-Schwarz inequality to obtain
\eq{\fr{2}{F^{2}}g^{ij}\ip{x_{;i}}{e}F_{;j}-\fr{2\ip{x}{e}}{F^{3}}\|\n F\|^{2}\leq \const.}
Hence, due to the concavity of $F$, the bounds on $F$ and the lower bound on $\ip{x}{e}$, the term involving $\|\cW\|^{2}$ is the leading order term with a negative sign. Due to the boundary condition \eqref{BD-H}, we obtain the result.
}

We obtain the long-time existence of the solution to \eqref{Flow}.

\begin{thm}\label{thm3.1}
Let 
\eq{x_{0}\cn\bar{\B}^{n}\hra\Si\sub\bar{\mathbb{B}}^{n+1}} be the embedding of a strictly convex free boundary hypersurface in the unit ball.
Then the maximal solution of \eqref{Flow} exists for all times with uniform $C^{\8}$-estimates.
\end{thm}

\pf{
In order to apply parabolic regularity theory, it is convenient to transform the flow to a scalar parabolic Neumann problem. This was already completely performed in \cite{LambertScheuer:04/2016}. Our previous a priori estimates allow us to use the same conclusion is in that paper. For convenience we sketch the main ingredients again. We use a coordinate transformation, which we call {\it{Moebius coordinates}} for the upper half-ball. Namely, consider the map
\eq{\p\cn \mathbb{B}^{n}\x [1,\8)&\ra \bar{\mathbb{B}}^{n+1}_{+}\\
				\p(x,\la)&=\fr{4\la x+(1+|x|^{2})(\la^{2}-1)e}{(1+\la^{2})+(1-\la)^{2}|x|^{2}}.}
This is a conformal map from the Euclidean cylinder $\bar{\mathbb{B}}^{n}\x[1,\8)$ to $\bar{\mathbb{B}}^{n+1}_{+}$ and in \cite[Prop~1]{LambertScheuer:04/2016} it was shown that any convex hypersurface of $\bar{\mathbb{B}}^{n+1}$ with perpendicular Neumann condition can be written as a graph over $\bar{\mathbb{B}}^{n}$ in Moebius coordinates with uniform gradient estimates, as long as $\la$ stays away from infinity. Hence all flow hypersurfaces can be written as graphs
\eq{\Si_{t}=\{(x(t,z),u(t,x(t,z)))\cn (t,z)\in [0,T^{*})\x \mathbb{B}^{n}\},}
where $u$ solves the Neumann problem
\eq{\dot{u}=-\fr{v}{e^{\psi}F}\q&\text{in}~(0,T^{*})\x\mathbb{B}^{n}\\
		\ip{Du}{N}=0\q&\text{on}~[0,T^{*})\x\del\mathbb{B}^{n}\\
			u(0,\cdot)=u_{0}\q&\text{on}~\{0\}\x\mathbb{B}^{n},}
where $v^{2}=1+|Du|^{2}$, $N$ is the outward normal to $\mathbb{B}^{n}$ and $e^{2\psi}$ is the conformal factor arising from the Moebius coordinates, compare \cite[Cor.~4]{LambertScheuer:04/2016}.
As in \cite[Lemma~15]{LambertScheuer:04/2016}, we can conclude the uniform $C^{\8}$-estimates and the long-time existence from \cite[Thm.~14.23]{Lieberman:/1998} or \cite[Thm.~4, Thm.~5]{Uraltseva:06/1994}. Note that we have also used that due to $F\geq c>0$ and $H\leq c$ the operator is uniformly parabolic.
}	

We finish the convergence result by proving that any subsequential limit is a spherical cap of a uniquely determined curvature.

\begin{prop}
The flow \eqref{Flow} converges smoothly to a uniquely determined spherical cap around $e$.
\end{prop}

\pf{\eq{\del_{t}W_{0}(\wh{\Si_{t}})=\int_{\Si_{t}}f\geq \int_{\Si_{t}}\br{\fr{n\ip{x}{e}}{H}-\ip{X}{\nu}}\geq 0.}
Furthermore
\eq{W_{0}(\wh{\Si_{t}})=|\wh{\Si_{t}}|}
is bounded and thus due to the $C^{\8}$-estimates we must have
\eq{\int_{\Si_{t}}\br{\fr{n\ip{x}{e}}{H}-\ip{X}{\nu}}\ra 0}
as $t\ra \8$.
Hence any limit satisfies the Heintze-Karcher type inequality with equality and thus must be a spherical cap due to \cite[Thm.~5.2]{WangXia:08/2017}. 
Due to the barrier property this cap is uniquely determined.
}

\section{Proof of the main results}\label{Main-Proof}
\subsection*{Proof of {\texorpdfstring{\cref{Main}}{}}}
(i)~First we prove that the functions
\eq{\label{Cap-Flow}f_{k}(r)=W_{k}(\wh{C_{r}(e)})}
are strictly increasing. This can be seen from the flow
\eq{\dot{y}&=X_{e}\\
	y(0,\bar \B^{n})&=C_{r}(e).}
The flow hypersurfaces of this flow are $(C_{s})_{s\geq r}$, where $s=r(t)$ is some increasing function and there holds
\eq{(n+1)\binom{n}{k-1}\del_{t}W_{k}(\wh{C_{r(t)}})=\int_{C_{r(t)}}k\s_{k}\ip{X_{e}}{\nu}>0}
due to the strict convexity of the caps. 

(ii)~We prove \cref{Main} in the strictly convex case.
Start the flow \eqref{Flow} from $\Si$ and denote the limit cap by $C_{R_{0}}$. Due to \cref{Monotone} there hold
\eq{W_{n}(\wh\Si)=W_{n}(\wh{C_{R_{0}}})=f_{n}(R_{0})=f_{n}\circ f_{k}^{-1}(f_{k}(R_{0}))\geq f_{n}\circ f_{k}^{-1}(W_{k}(\wh\Si)).} This inequality is strict, unless $\Si$ is a spherical cap, due to the final statement in \cref{Monotone}.

(iii)~We consider the convex case, where we may suppose that $\Si$ is not a flat disk, otherwise the statement is trivially true.
By \cite[Cor.~3.3]{LambertScheuer:09/2017} we can approximate $\Si$ by strictly convex hypersurfaces with free boundary in the unit ball in the $C^{2,\al}$-norm. Hence the inequality holds for $\Si$.

To prove the limiting case, we employ an argument previously used in \cite{GuanLi:08/2009}. Suppose $\Si$ is convex and
\eq{W_{n}(\wh\Si)=f_{n}\circ f_{k}^{-1}(W_{k}(\wh\Si)).}
Due to \cite[Lemma~3.1]{LambertScheuer:09/2017}
the set
\eq{\Si_{+}=\{\xi\in \B^{n}\cn (h_{ij}(\xi))>0\}}
is not empty, since $\Si$ must have a strictly convex point in the interior. $\Si_{+}$ is obviously open. We prove that $\Si_{+}$ is also closed, by showing that 
\eq{h_{ij|\Si_{+}}\geq cg_{ij|\Si_{+}}}
for some constant $c$.	
So let $\xi_{0}\in \Si_{+}$, $U\sub\B^{n}$ open with
\eq{\xi_{0}\in U\sub\bar U\sub \B^{n},}
 $\eta\in C^{\8}_{c}(U)$ and consider the normal variation
\eq{x_{s}=x_{0}+s\eta\nu,}
where $x_{0}$ is the embedding of $\Si$. Obviously, for small $s$ the corresponding hypersurfaces $\Si_{s}$ are convex with free boundary in the unit ball and hence
\eq{W_{n}(\wh{\Si_{s}})\geq f_{n}\circ f_{k}^{-1}(W_{k}(\wh{\Si_{s}})).}
Thus
\eq{0&=n(n+1)\fr{d}{ds}_{|s=0}\br{W_{n}(\wh{\Si_{s}})-f_{n}\circ f_{k}^{-1}(W_{k}(\wh{\Si_{s}}))}=\int_{\Si}\br{n\s_{n}-c_{1}\s_{k}}\eta}
with a suitable positive $c_{1}$.
Since $\eta$ was arbitrary, we have on all of $\B^{n}$:
\eq{n\s_{n}=c_{1}\s_{k},}
where
Since 
\eq{\s_{n}\leq c_{n,k}\s_{k}^{\fr nk},}
we obtain
\eq{\s_{n}=\fr{c_{1}}{n}\s_{k}\geq \fr{c_{1}}{n}\br{\fr{\s_{n}}{c_{n,k}}}^{\fr kn}.}
Hence $\s_{n}$ has a positive lower bound on $\Si_{+}$, which implies that also the second fundamental form has a positive lower bound on $\Si_{+}$. Thus $\Si_{+}$ is closed. Hence $\Si$ is strictly convex and by the first case it is a spherical cap.

\section{Appendix A. Quermassintegrals and the Euler characteristic}

In this Appendix, we show that $W_{n+1}$ is a topological constant for a general free boundary hypersurface in the unit ball.
First  we can easily extend the definition of $W_{n+1}$ as follows:
\enum{
\item
If $n$ is even and $\Si$ immersed, the definition of quermassintegrals given in \eqref{Def-QM} carries over, since it only depends on the area of $\Sigma$ and curvature  integrals.
We still denote it by $W_{n+1}(\wh{\Si})$, though
$\wh{\Sigma}$ is not well-defined. 
\item
If $n$ is odd, we consider an embedded hypersurface $\Sigma\subset \bar\B^{n+1}$ with an embedded free boundary $\partial \Sigma \subset \bbS^{n}$. 
Since $\partial \Sigma  \subset \bbS^{n}$
is an embedding, it divides the sphere into 2 pieces. Let us choose one of these 2 pieces and denote it as before by $\wh{\del\Si}$. Then
$\wh{\Si}$ is the domain in $\bar\B^{n+1}$ enclosed by $\wh{\del\Si}$ and $\Sigma$ and $W_{n+1}(\wh\Si)$ is defined as in \eqref{Def-QM}.\footnote{Note that if we choose the other piece, we have to flip the normal in the definition of the principal curvatures.}
}

Moreover, the first variation formula of $W_{n+1}$ is the same as the one proved Proposition \ref{Lemma-Wk}, namely

\eq{
 \partial_t W_{n+1}( \wh{\Sigma})=0
}
for any smooth deformation.
This formula implies the following Gauss-Bonnet-Chern theorem.
\begin{thm}\label{GB_thm}
 Let $\Sigma\subset \bar\B^{n+1}$ be  immersed hypersurface with  free boundary $\partial \Sigma \subset \bbS^{n}$.
 If $n$ is even, then
 \eq{ 
   \frac {\omega_n} 2 \chi(\Sigma)&= (n+1)W_{n+1}(\wh\Si)\\
      &= \int_\Sigma H_n+ \int_{\partial\Si}H_{n-2}^{\bbS^n}
      +\frac {n-2} 3 \int_{\partial\Si}H_{n-4}^{\bbS^n} +\cdots \\
      &\hp{=} +\frac{n-2} 3 \cdots \frac 4{n-3}  \int _{\partial\Si}H_{2}^{\bbS^n}    + \frac{n-2} 3 \cdots \frac 2{n-1} |\partial \Si|.
     }
 If $n$ is odd and $\Sigma\subset \bar\B^{n+1}$ is embedded with an embedded boundary $\partial \Si$, then
 \eq{
  {\omega_n} \chi(\wh{\Sigma})&=  (n+1) W_{n+1}(\wh{\Sigma}) \\
  &= \int_\Sigma H_n+\int_{\partial\Si}H_{n-2}^{\bbS^n}
      \hp{=}+\frac {n-2} 3 \int_{\partial\Si}H_{n-4}^{\bbS^n} +\cdots + \frac{n-2}{3}   \int _{\partial\Si}H_{1}^{\bbS^n}    +  |\wh {\partial \Si}|.
  }
   Here $H_k$ is the (normalized) $k$-th mean curvature of $\Sigma$ in $\bbR^{n+1}$
   and $H_{j}^{\bbS^n}$ is the (normalized) $j$-th mean curvature of $ \partial \Sigma$ in $\bbS^n$.
\end{thm}

 Note that in the second case $\chi(\wh{\Sigma})=  \frac 12 \chi(\Si\cup \wh{\partial\Si})$.

\pf{The second equality follows from the very definition. We need to prove the first, which follows from the variation formula \eqref{Wk-top}
of $W_{n+1}$ by a standard argument as follows.
One deforms the hypersurface $\Si$ to a very small free boundary hypersurface $\widetilde \Si$. \eqref{Wk-top} implies that both $\Sigma$ and $\widetilde \Si$
have the same $W_{n+1}$. It is easy to see that except the leading order term, all other terms of $W_{n+1}$ become small, 
when $\widetilde \Si$ shrinks to a point. In this case,  $\widetilde \Si$ looks more and more like a free boundary hypersurface 
in the halfspace $\bbR^{n+1}_+$ with the boundary perpendicular on $\partial \bbR^{n+1}_+$. Hence the first equality follows from the Gauss-Bonnet theorem for 
closed hypersurfaces in $\bbR^{n+1}.$
}

We would like to compare the above Gauss-Bonnet-Chern theorem to the Gauss-Bonnet theorem for closed hypersurfaces in $\bbR^{n+1}$ by Hopf \cite{Hopf:/1926,Hopf:/1927a} and in $\bbS^n$ by Teufel \cite{Teufel:/1980}. For the Gauss-Bonnet theorem for closed hypersurfaces in the hyperbolic space, we refer to \cite{Solanes:/2006}.

If  $n$ is even and $\Si$ is an immersed closed hypersurface in $\bbR^{n+1}$, then by \cite{Hopf:/1926}, we have 
\[  \frac {\omega_n} 2 \chi{(\Si)}=(n+1) W^{\bbR^{n+1}}_{n+1}=\int_{\Sigma} H_n .\]
If $n$ is odd and $\Si$ is an embedded closed hypersurface in $\bbR^{n+1}$, then by \cite{Hopf:/1927a}, we have 
\[ {\omega_n} \chi{(\Omega)}=(n+1) W^{\bbR^{n+1}}_{n+1}=\int_{\Si} H_n,\]
where $\Omega$ is the domain enclosed by $\Si$. Note that in this case $\chi{(\Omega)}=\frac 12 \chi(\Si)$.


For  a closed hypersurfaces in $\bbS^n$,  Teufel \cite{Teufel:/1980} obtained: If $n$ is odd and  $\Sigma$ is immersed,
then 
\[
 \frac {\omega_{n-1}} 2 \chi{(\Si)}= 
 n W^{\bbS^n}_n (\Sigma)= c_{n-1}\int_\Sigma H_{n-1} +c_{n-3}\int_\Sigma H_{n-3}+\cdots +  c_2 \int_\Sigma H_2 + |\partial \Si|.
\]
If $n$ is even and  $\Sigma$ is embedded with the enclosed domain $K$, then 
\[
 \om_{n-1}\chi{(K)}= 
 n W^{\bbS^n}_n (K)= c_{n-1}\int_\Sigma H_{n-1} +c_{n-3}\int_\Sigma H_{n-3}+\cdots +  c_1 \int_\Sigma H_1 + |K|.
\]
For the constants $c_j$, we refer to  \cite{Teufel:/1980}.


\subsection*{Acknowledgements}
Large parts of this work were produced during a research visit of JS to the school of mathematical sciences at Xiamen university. JS would like to express his deep gratitude to the math department in Xiamen and especially to CX for their hospitality and generosity.

\noindent
CX is supported by NSFC (Grant No. 11871406), the Natural Science Foundation of Fujian Province of China (Grant No. 2017J06003) and the Fundamental Research Funds for the Central Universities (Grant No. 20720180009).

\bibliographystyle{hamsplain}
\bibliography{Bibliography}

\end{document}